\documentclass[a4paper,
               headinclude=false,
               headsepline,
               footinclude=false,
               footsepline,
               plainfootsepline,
               abstract=true,
               twoside]{scrartcl}
\usepackage[active]{srcltx} 
\usepackage{scrlayer-scrpage}
\RequirePackage{hyperref}
\usepackage{amsmath,amssymb,amsthm}
\usepackage{dsfont}
\usepackage{bbm}
\usepackage{libertinus}\providecommand{\coloneqq}{:=}
\usepackage[svgnames]{xcolor}
\usepackage{tikz}
\usetikzlibrary{shapes}
\usetikzlibrary{calc}
\addtokomafont{title}{\rmfamily}
\addtokomafont{section}{\rmfamily}
\addtokomafont{subsection}{\rmfamily}
\swapnumbers%% and make sure that numbers are swapped by not using
\theoremstyle{plain}
\newtheorem{theo}{Theorem}[section]
\newtheorem{coro}[theo]{Corollary}
\newtheorem{lemm}[theo]{Lemma}
\newtheorem*{lemm*}{Lemma}
\newtheorem{prop}[theo]{Proposition}
\newtheorem{namet}[theo]{\myThmName}
\newtheorem*{namet*}{\myThmName}

\newenvironment*{nthm*}[1][\kern-.35em]{\def\myThmName{#1}\begin{namet*}}{\end{namet*}}

\theoremstyle{definition}
\newtheorem{defi}[theo]{Definition}
\newtheorem*{defi*}{Definition}
\newtheorem{nota}[theo]{Notation}

\newtheorem{exas}[theo]{Examples}
\newtheorem{rema}[theo]{Remark}
\newtheorem*{rema*}{Remark}
\newtheorem{rems}[theo]{Remarks}
\newtheorem{named}[theo]{\myThmName}
\newtheorem*{named*}{\myThmName}
\newenvironment{ndef}[1][\kern-.35em]{\def\myThmName{#1}\begin{named}}{\end{named}}

\newcounter{rememberEnumi}
\newcommand{\SaveEnumi}{\global\setcounter{rememberEnumi}{\value{enumi}}}
\newcommand{\RecallEnumi}{\setcounter{enumi}{\therememberEnumi}}
\newcommand{\set}[2]{\left\{{#1}\left|\vphantom{#1#2\strut}\right.\, 
                    {#2}\right\}}

\newcommand{\SL}[2]{\mathrm{SL}(#1,#2)}
\newcommand{\GL}[2]{\mathrm{GL}(#1,#2)}
\newcommand{\gL}[2]{\mathrm{\Gamma L}(#1,#2)}
\newcommand{\sL}[2]{\mathrm{\Sigma L}(#1,#2)}
\newcommand{\PSL}[2]{\mathrm{PSL}(#1,#2)}
\newcommand{\PGL}[2]{\mathrm{PGL}(#1,#2)}
\newcommand{\PgL}[2]{\mathrm{P\Gamma L}(#1,#2)}

\newcommand{\Gtwo}{\mathrm{G}_2}
\newcommand{\tGtwo}{\mathrm{^2G}_2}

\newcommand{\gU}[2]{\mathrm{\Gamma U}(#1,#2)}
\newcommand{\PgU}[2]{\mathrm{P\Gamma U}(#1,#2)}

\newcommand{\Sym}[2][]{\mathrm{Sym}(#2)}
\newcommand{\Aut}[2][]{\mathrm{Aut}_{#1}(#2)}
\newcommand{\Fix}[2][]{\mathrm{Fix}_{#1}(#2)}
\newcommand{\C}[2][]{\mathrm{C}_{#1}(#2)}%% centralizer
\newcommand{\N}[2][]{\mathrm{N}_{#1}(#2)}%% normalizer
\newcommand{\p}[1]{\hat{#1}}

\newcommand{\mybb}{\mathds}
\newcommand{\1}{\mybb1}
\newcommand{\Eta}{\mathrm{H}}
\newcommand{\Zeta}{\mathrm{Z}}

\newcommand{\FF}{\mybb F}
\newcommand{\KK}{\mybb K}
\newcommand{\Ss}{\mybb S}

\newcommand{\Char}{\operatorname{char}}

\newcommand{\id}{\operatorname{id}}
\newcommand{\tr}{\operatorname{tr}}

\newcommand{\cB}{\mathcal B}
\newcommand{\cJ}{\mathcal J}
\newcommand{\cK}{\mathcal K}
\newcommand{\cN}{\mathcal N}
\newcommand{\cS}{\mathcal S}

\newcommand{\elt}[3]{\left(#1,#2,#3\right)^\intercal}
\newcommand{\Elt}[4][0ex]{{\left(\begin{array}{c}#2\\[#1]#3\\[#1]#4\end{array}\right)}}%
\newcommand{\centeredEntry}[3]{{}#1{}&{}#2{}&{}#3{}}
\newcommand{\EltAl}[4][0ex]{{\arraycolsep=0pt\left(\begin{array}{rcl}
                                      \centeredEntry#2\\[#1]
                                      \centeredEntry#3\\[#1]
                                      \centeredEntry#4
                                    \end{array}\right)}}%
\newcommand{\RU}[2][]{\Xi(\ifx\empty#1\else#1,\fi#2)}
\newcommand{\RUold}[2][]{\mathrm{U}(\ifx\empty#1\else#1,\fi#2)}
\newcommand{\RT}[2][]{\mathrm{RT}(\ifx\empty#1\else#1,\fi#2)}
\newcommand{\Ree}[2][]{\mathrm{Ree}(\ifx\empty#1\else#1,\fi#2)}

\DeclareFontFamily{U}{mathb}{\hyphenchar\font45}%
\DeclareFontShape{U}{mathb}{m}{n}{%
<-6> mathb5 <6-7> mathb6 <7-8> mathb7 %
<8-9> mathb8 <9-10> mathb9 %
<10-12> mathb10 <12-> mathb12 %
}{}%
\DeclareSymbolFont{mathb}{U}{mathb}{m}{n}%
\DeclareMathSymbol{\llcurly}{\mathrel}{mathb}{"CE}%
  \newcount\Hour
   \newcount\Minute
   \def\timenow{\Hour=\time \Minute=\Hour 
        \divide\Hour by 60 \number\Hour:%
        \multiply\Hour by 60 %\@minute=\time 
        \global\advance\Minute by -\Hour%
        \ifnum\Minute<10 0\number\Minute\else\number\Minute\fi}
\let\tOday\today \let\today\empty
\lohead[]{\footnotesize\normalfont\rmfamily \MYtitle}
\lehead[]{\footnotesize\normalfont\rmfamily \MYauthor}
\rehead[]{\footnotesize\normalfont\rmfamily \MYtitle}
\rohead[]{\footnotesize\normalfont\rmfamily \MYauthor}
\lefoot[\pagemark]{\pagemark}
\cefoot[\footnotesize\normalfont\ttfamily\tiny version of \tOday]%
       {\footnotesize\normalfont\ttfamily\tiny version of \tOday}
\refoot[\tiny\copyright~{\normalfont\ttfamily
  by~\MYauthor
}]{\tiny\copyright~{\normalfont\ttfamily by~\MYauthor}}
\lofoot[\tiny\copyright~{\normalfont\ttfamily
  by~\MYauthor
}]{\tiny\copyright~{\normalfont\ttfamily by~\MYauthor}} 
\cofoot[\footnotesize\normalfont\ttfamily\tiny version of \tOday]%
       {\footnotesize\normalfont\ttfamily\tiny version of \tOday}
\rofoot[\pagemark]{\pagemark}
\newcommand{\keywords}[1]{\par\noindent{\normalfont\bfseries Keywords: }#1}
\newcommand{\subjclass}[1]{\par\noindent{\normalfont\bfseries Mathematics Subject
    Classification (MSC 2000): }#1}
\newcommand{\MSC}[1]{\href{https://mathscinet.ams.org/mathscinet/freetools/msc-search?text=#1}{#1}}
\title{Construction of the smallest Ree-Tits unital\\
  from the special linear group of degree two over the field with
  eight elements}%
\author{Markus J. Stroppel}
  \let\MYauthor\shortauthor 
  \let\MYtitle\shorttitle
\begin{document}
\maketitle
%%%%%%%%%%%%%%%%%%%%%%%%%%%%%%%%%%%%%%%%%%%%%%%%%%%%%%%%%%%%%%%%%%%%%%
\begin{abstract}\noindent %
  We construct the smallest Ree-Tits unital from a group of matrices that
  is isomorphic to the commutator group of the corresponding Ree
  group.
  The matrix description is used to determine configurations in the
  unital via explicit computations. %
  Embeddings into larger Ree-Tits unitals are made explicit. 
\end{abstract}
%%%%%%%%%%%%%%%%%%%%%%%%%%%%%%%%%%%%%%%%%%%%%%%%%%%%%%%%%%%%%%%%%%%%%%
\subjclass{%
  \MSC{51E26}, % Other finite linear geometries
  \MSC{20D06}, % Simple groups: alternating groups and groups of Lie type
}%
\keywords{Ree group, classical group, exceptional isomorphism,
  unital, automorphism, O'Nan configuration, Ree-Tits unital }
%%%%%%%%%%%%%%%%%%%%%%%%%%%%%%%%%%%%%%%%%%%%%%%%%%%%%%%%%%%%%%%%%%%%%%

\section*{Introduction}

As one of the last infinite families of finite simple groups, the Ree
groups of type $\Gtwo$ were introduced in~\cite{MR0138680}. See
also~\cite{MR1611778} and~\cite{MR0193136}. In part of the literature,
these groups are denoted by~$\tGtwo(3^{2\ell+1})$; they are simple
groups for each positive integer~$\ell$, but $\tGtwo(3)$ has a normal
subgroup of index~$2$. That normal subgroup is isomorphic to~$\SL28$.

Let~$\ell$ be a non-negative integer, and put
$q \coloneqq 3^{2\ell+1}$.  There is a combinatorial geometry
associated with the Ree group~$\tGtwo(q)$, namely, a unital of
order~$q$ (viz., a $2${-}$(q^3+1,q+1,1)$ design) known as the
\emph{Ree-Tits unital}, see~\cite{MR0193136}. %
Using the classification of finite simple groups,
Kantor~\cite{MR773556} has shown that the Hermitian unitals
(\cite[p.\,104]{MR0233275}, %
\cite[2.1, 2.2, see also p.\,29]{MR2440325}) and the Ree-Tits unitals
are the only unitals that admit a group of automorphisms that is doubly
transitive on the set of points.  For the special case of unitals of
order~$3$, this was found in a computer search by
Brouwer~\cite{MR655065}.  Many isomorphism classes of unitals of
order~$3$ are known, see~\cite{MR1991559}. %

In the present paper, we give a simple direct construction of the
smallest Ree-Tits unital (namely, the one of order~$3$) in terms of the
group $\SL28$. We use this description to compute explicitly some
joining blocks and intersection points. In particular, we exhibit an
O'Nan configuration, and actually, a super O'Nan configuration (viz.,
the dual of the complete graph on~$5$ vertices).
We translate the findings into general Ree-Tits unitals. %
For the readers' convenience, we include the relevant information on
Ree groups and Ree-Tits unitals. 

\section{The special linear group}
\label{sec:theGroup}

\begin{nota}\label{nota:PgL28}
  The field $\FF_8$ of order~$8$ is obtained by adjoining a root~$u$
  of $X^3+X+1$ to~$\FF_2$. The elements of~$\FF_8$ are then $0$, $1$,
  $u$, $u^2$, $u^3=u+1$, $u^4=u^2+u$, $u^5=u^2+u+1$, and $u^6=u^2+1$.
  Note that $u^7=1$; each element of $\FF_8\smallsetminus\{0\}$ is a root
  of $X^7-1 = ({X-1})({X^3+X+1})({X^3+X^2+1})$.  The roots of $X^3+X+1$ are
  $u$, $u^2$, and $u^4$.  The roots of $X^3+X^2+1$ are $u^3$, $u^6$,
  and $u^5 = (u^3)^4$.

  The group $\gL28$ of all semi-linear bijections is the semi-direct
  product $\Aut{\FF_8}\ltimes\GL28$. %
  Let $\delta\in\gL28$ be defined by
  $(x_0,x_1)^\delta \coloneqq (x_0^4,x_1^4)$, then conjugation
  by~$\delta$ maps $x = \left(
    \begin{smallmatrix}
      x_{00} & x_{01} \\
      x_{10} & x_{11} 
    \end{smallmatrix}\right) \in \GL28$ to 
  \[
    x^\delta \coloneqq \delta^{-1}x\delta = %
    \left(
      \begin{matrix}
        x_{00}^4 & x_{01}^4 \\[1ex]
        x_{10}^4 & x_{11}^4 
      \end{matrix}\right).
  \]
  We find $\gL28 = \langle\delta\rangle\GL28$. %
  In $\SL28$, consider
  \[
    \1 \coloneqq \left(
      \begin{matrix}
        1 & 0 \\
        0 & 1 
      \end{matrix}\right),
    \quad S \coloneqq \left(
      \begin{matrix}
        0 & 1 \\
        1 & 0 
      \end{matrix}\right),
    \quad T \coloneqq \left(
      \begin{matrix}
        1 & 0 \\
        1 & 1 
      \end{matrix}\right),
    \quad A \coloneqq \left(
      \begin{matrix}
        u^2 & u \\
        u & u^4
      \end{matrix}\right),
    \quad D \coloneqq \left(
      \begin{matrix}
        0 & 1 \\
        1 & 1
      \end{matrix}\right).
  \]
  Then $A^3 = D = TS$, and $A^4 = A^\delta$.  For $M\in\SL28$, we note
  that $M^S = S^{-1}MS$ is the transpose of the inverse of~$M$. %
  In particular, we have $A^S = A^{-1}$.
\end{nota}

\begin{lemm}
  An element\/ $M\in\SL28\smallsetminus\{\1\}$ has
  \begin{itemize}
  \item order~$2$ if\/ $\tr(M) = 0$,
  \item order $3$ if\/ $\tr(M)=1$,
  \item order $7$ if\/ $\tr(M)\in\{u^3,u^6,u^5\}$,
  \item order $9$ if\/ $\tr(M)\in\{u,u^2,u^4\}$.
  \end{itemize}
  In particular, the order of\/~$A$ is\/~$9$. %
\end{lemm}
\begin{proof}
  The center of $\SL28$ is trivial because $s^2=1$ implies $s=1$ since
  the ground field has characteristic~$2$.  As all elements in
  question have determinant~$1$, the trace~$\tr(M)$ determines~$M$, up
  to conjugacy in~$\GL28$, and elements with traces in the
  $\Aut{\FF_8}$-orbits~$\{u,u^2,u^4\}$ or $\{u^3,u^6,u^5\}$ have
  conjugates in~$\gL28$ with trace~$u$ or~$u^3$, respectively.  Now it
  suffices to check the claim for the representatives $S$, $D$, $A$
  and the element $\left(
    \begin{smallmatrix}
      u & 0 \\
      0 &u^6 
    \end{smallmatrix}\right)$ with trace~$u^5$ and order~$7$. 
\end{proof}

\begin{rems}\label{sylowPsL28}
  The orders of the pertinent groups are
  $|\GL28| = (8^2-1)(8^2-8) = 2^3\cdot3^2\cdot7^2$,
  $|\SL28| = |\PGL28| = |\GL28|/7 = 2^3\cdot3^2\cdot7$, and
  $|\PgL28| = 3\cdot|\PGL28| = 2^3\cdot3^3\cdot7$. %
  As the center of $\SL28$ is trivial, we obtain
  $\SL28 \cong \PSL28 \cong \PGL28$, and
  $\sL28 \coloneqq \langle\delta\rangle\SL28 \cong \PgL28$. %
  (E.g., see the first chapter of~\cite{MR1859189} for these well
  known facts). %

  A Sylow $3$-subgroup~$\Delta$ of $\gL28$ is generated by $\delta$
  and~$A$. That group is a semidirect product
  $\langle\delta\rangle\ltimes\langle A\rangle$; we already know
  $A^\delta = A^4$. For $m\in\{0,1,2\}$, we note
  $4^{2m}+4^m+1 \equiv 3 \pmod9$, and conclude
  $(\delta^mA^n)^3 = \delta^{3m}A^{n(4^{2m}+4^m+1} = A^3$. Therefore,
  the elements of order~$9$ in~$\Delta$ are the~$18$ elements of the
  form $\delta^m A^n$ where $n\in\{1,2,4,5,7,8\}$ and $m\in\{0,1,2\}$;
  every other non-trivial element of~$\Delta$ has order~$3$. The
  subgroup $\langle\delta,A^3\rangle$ is elementary abelian. The
  center of~$\Delta$ is~$\langle A^3\rangle$, and consists of all
  third powers of elements of~$\Delta$.  Note also that
  $\SL28\cap \Delta = \langle A\rangle$, that
  $\Delta = \C[\sL28]{A^3}$, and that
  $\p{A} \coloneqq \langle S\rangle\ltimes\Delta =
  \langle\{S\}\cup\Delta\rangle = \N[\sL28]\Delta = \N[\sL28]{\langle
    A^3\rangle}$. %
  We infer that $\sL28$ has $|\sL28|/(2|\Delta|) = %8\cdot9\cdot7 =
  4\cdot7 = 28$ Sylow $3$-subgroups.

  Any two given Sylow $3$-subgroups of~$\sL28$ have trivial
  intersection because their sets of third powers are
  different. Thus~$\Delta$ acts semi-regularly, and then sharply
  transitively on the set of all other Sylow $3$-subgroups,
  and~$\sL28$ is doubly transitive on the set of its Sylow
  $3$-subgroups.

  As $\sL28/\SL28$ is a group of order~$3$, every involution in
  $\sL28$ actually lies in $\SL28$, and is a conjugate of~$S$ in~$\SL28$. The
  centralizer of~$S$ in~$\SL28$ is $\C[\SL28]S =
  \1 +\FF_8\left(
    \begin{smallmatrix}
      1 & 1 \\
      1 & 1 
    \end{smallmatrix}\right)$.
  We obtain that $\sL28$ contains exactly $|\SL28|/8 = 63$
  involutions. %
  The centralizer of~$S$ in $\N[\sL28]\Delta$ is
  $\langle S,\delta\rangle$, and we find that $\langle A\rangle$ acts
  transitively on the set of~$9$ involutions in~$\N[\sL28]\Delta$.
\end{rems}

\section{The unital}

\begin{defi}
  We form the incidence geometry $\Ss \coloneqq (\cN,\cJ,\ni)$,
  where~$\cN$ is the set of all normalizers of Sylow $3$-subgroups,
  and~$\cJ$ is the set of all involutions in~$\sL28$.  %
  Recall that the normalizer of any Sylow $3$-subgroup in~$\sL28$
  equals the normalizer of the center of that subgroup, and also
  equals the normalizer of the intersection of that subgroup
  with~$\SL28$.  %
  (In Section~\ref{sec:isoPsL28Ree3} below, we will see that~$\Ss$ is
  isomorphic to the Ree-Tits unital~$\RT3$ of order~$3$.) %
\end{defi}

\begin{theo}
  The incidence geometry\/~$\Ss$ is a unital of order~$3$.  The group
  $\SL28$ acts flag-transitively on~$\Ss$.
\end{theo}
\begin{proof}
  Consider $\p{A} \coloneqq \N[\sL28]\Delta \in\cN$. %
  From~\ref{sylowPsL28} we know that~$\SL28$ acts transitively both
  on~$\cN$ and on~$\cJ$, and acts by automorphisms of~$\Ss$. %
  We have also seen that $|\cN| = 28 = 3^3+1$, $|\cJ|=63$, and
  $|N\cap\cJ| = 9$. %
  So every point is incident with exactly~$9$ blocks, and there are
  $28\cdot 9$ flags. We obtain the number of points per block as
  $(28\cdot9)/63 = 4$.

  If $(N^g,S^h)$ is a flag then $S^{hg^{-1}}$ is an involution in~$N$,
  and there exists $k\in N$ with $S^k = S^{hg^{-1}}$. Now $N^k = N$,
  and $(N,S)^{kg} = (N^{kg},S^{kg}) = (N^g,S^h)$, as required.

  It remains to show that any two points are on a unique block.
  Consider two points $N^g$ and $N^h$ in~$\cN$. Without loss of
  generality, we may assume $N^h =N$.
  Each non-trivial element of $N' = \langle A\rangle$ generates the
  subalgebra $\FF_8[A] \cong \FF_{64}$ in the endomorphism ring
  of~$\FF_8^2$. That subalgebra has dimension~$2$ over~$\FF_8$, and
  $A^g$ generates a different subalgebra. So the set
  $({N^g\cap N})\smallsetminus\{\1\}$ is contained in $N\smallsetminus N'$, and
  consists of involutions. As~$N$ is a dihedral group of order
  $2\cdot9$, we obtain $|N^g\cap N|\le2$.  This means that there is
  at most one block joining~$N$ and~$N^g$. %
  There are~$9$ blocks through~$N$, each one contains~$3$ points apart
  from~$N$, and none of those is on two blocks through~$N$. Thus
  $1+9\cdot 3 = 28$ points are joined to~$N$: these are all the
  points, as required.
\end{proof}

Joining blocks and intersections in~$\Ss$ are fairly easy to compute:
\begin{lemm}\label{joinAndMeet}\quad
  \begin{enumerate}
  \item Let\/ $B,C\in\SL28$ be elements of order~$3$ or~$9$ such
    that\/ $\N[\SL28]{\langle B\rangle} \ne %
    \N[\SL28]{\langle C\rangle}$. Then these two points of\/~$\Ss$ are
    incident with a unique common block~$I\in\cJ$, namely, the unique
    involution $I \in %
    \N[\SL28]{\langle B\rangle} \cap %
    \N[\SL28]{\langle C\rangle}$. %
    That involution can be found by solving the two conditions
    $IB = B^{-1}I$ and $IC = C^{-1}I$ simultaneously.
  \item For $I,L\in\cJ$, there exists a point incident with both if,
    and only if, the product $IL$ has order\/~$3$ or\/~$9$, viz., if\/
    $\tr(IL) \in \{1,u,u^2,u^4\}$. \\
    That point is then unique, it is the normalizer
    of\/~$\langle IL\rangle$. %
    \qed
  \end{enumerate}
\end{lemm}

\begin{theo}
  The unital\/~$\Ss$ is not isomorphic to a Hermitian unital. 
\end{theo}
\begin{proof}
  We give a group-theoretic argument; one could also use the fact
  that~$\Ss$ contains an O'Nan configuration
  (see~\ref{explicit}\ref{oNan} below) while the Hermitian unitals do
  not contain such configurations (see~\cite[Proposition,
  p.\,507]{MR0295934}, cp. also~\cite[2.2]{MR2795696}).

  In~\ref{sylowPsL28}, we have seen that the Sylow $3$-subgroups of
  $\SL28$ are cyclic of order~$9$. The full automorphism group of the
  Hermitian unital of order~$3$ is the group $\PgU3{\FF_9|\FF_3}$
  induced by the group $\gU3{\FF_9|\FF_3}$ of all semi-similitudes of
  a non-degenerate Hermitian form on~$\FF_9^3$ (\cite{MR0295934},
  see~\cite{MR2241352}). The Sylow $3$-subgroups of
  $\gU3{\FF_9|\FF_3}$ have exponent~$3$ (in fact, they are contained,
  up to conjugation, in groups of strictly upper triangular $3\times3$
  matrices, and the latter groups are Heisenberg groups,
  cp.~\cite[6.1, 2.3]{MR2410562}), so there are no elements of
  order~$9$ in $\gU3{\FF_9|\FF_3}$. %
  This implies that $\SL28$ is not contained in $\PgU3{\FF_9|\FF_3}$,
  and the unital $\Ss$ is not isomorphic to the Hermitian unital.
\end{proof}

\begin{rema}
  With arguments like those in~\cite[Section\,2]{MR847092} one can
  show that the unital~$\Ss$ is isomorphic to the %
  \emph{Ree-Tits unital} of order~$3$,
  compare~\cite{MR0193136}. Assuming the mere existence of an
  isomorphism from $\sL28$ onto $\Ree3$, we identify the unital in
  Section~\ref{sec:isoPsL28Ree3} below.
\end{rema}

\begin{rems}
  The group $\sL28$ acts two-transitively on~$\cN$, but $\SL28$ does
  not act two-transitively on~$\cN$; in fact, the order
  $|\cN\smallsetminus\{\p A\}| = 27$ does not divide the order
  $|\SL28\cap\p A| = 18$ of the stabilizer in~$\SL28$.
  The incidence geometry $\Ss$ can also be reconstructed as a coset
  geometry (\cite[6.2, 6.3]{MR797151}, \cite{MR0131216}); in
  particular, it is a sketched geometry for $\SL28$
  (see~\cite{MR1161931}, \cite{MR1214465}).
\end{rems}

%%%%%%%%%%%%%%%%%%%%%%%%%%%%%%%%%%%%%%%%%%%%%%%%%%%%%%%%%%%%%%%%%%%%%
\begin{exas}\label{explicit}
  We use the elements $\delta$, $S$, $T$, $A$, $D$ introduced
  in~\ref{nota:PgL28}, observe $S\delta = \delta S$, and abbreviate
  $\p X \coloneqq \N[\sL28]{\langle X\rangle}$.
  \begin{enumerate}
  \item %
    The blocks through~$\p D$ are the involutions $S$, $T$,
    $T^S = \left(
      \begin{smallmatrix}
        1 & 1 \\
        0 & 1 
      \end{smallmatrix}\right)$, %
    $Y \coloneqq \left(
      \begin{smallmatrix}
        u^{\phantom1} & u^2\\
        u^4 & u^{\phantom1} 
      \end{smallmatrix}\right)$, %
    $Y^\delta$, $Y^{\delta^2}$, %
    $Y^S = \left(
      \begin{smallmatrix}
        u^{\phantom1} & u^4\\
        u^2 & u^{\phantom1} 
      \end{smallmatrix}\right)$, %
    $Y^{S\delta}$, and~$Y^{S\delta^2}$. %
  \item Each one of the blocks $S$, $T$, and $T^S$ is fixed
    by~$\delta$, and so is the common point~$\p D$ on all of these
    blocks.
  \item The points on~$S$ are obtained as~$\p X$, where~$X$ is one of
    the following elements of order~$3$: \quad%
    $D$, $G \coloneqq \left(
      \begin{smallmatrix}
        u & u^6 \\
        u^6 & u^3
      \end{smallmatrix}\right)$, %
    $G^\delta = \left(
      \begin{smallmatrix}
        u^4 & u^3 \\
        u^3 & u^5
      \end{smallmatrix}\right)$, %
    $G^{\delta^2} = \left(
      \begin{smallmatrix}
        u^2 & u^5 \\
        u^5 & u^6
      \end{smallmatrix}\right)$. %
  \item The points on $T = \left(
      \begin{smallmatrix}
        1 & 0 \\
        1 & 1
      \end{smallmatrix}\right)$ are obtained as~$\p X$, where~$X$ is
    one of the following elements of order~$3$: \quad%
    $ D $, \ $ E \coloneqq \left(
      \begin{smallmatrix}
        1 & u^{\vphantom1} \\
        1 & u^3
      \end{smallmatrix}\right)^3
    = \left(
      \begin{smallmatrix}
        u^5 & 1 \\
        u^6 & u^4
      \end{smallmatrix}\right)$, \
    $E^\delta = \left(
      \begin{smallmatrix}
        u^6 & 1 \\
        u^3 & u^2
      \end{smallmatrix}\right)$, \
    and $ E^{\delta^2} = \left(
      \begin{smallmatrix}
        u^3 & 1 \\
        u^5 & u
      \end{smallmatrix}\right)$.
  \item The points on $T^S = \left(
      \begin{smallmatrix}
        1 & 1 \\
        0 & 1
      \end{smallmatrix}\right)$  are obtained as~$\p X$, where~$X$ is
    one of the following elements of order~$3$: \quad%
    $D$, \ $ E^S = \left(
      \begin{smallmatrix}
        u^4  & u^6 \\
        1 & u^5
      \end{smallmatrix}\right)$, \
    $E^{S\delta} = \left(
      \begin{smallmatrix}
        u^2 & u^3 \\
        1 & u^6
      \end{smallmatrix}\right)$, \
    and $E^{S\delta^2} = \left(
      \begin{smallmatrix}
        u & u^5 \\
        1 & u^3
      \end{smallmatrix}\right)$.
  \item %
    The block $M \coloneqq \left(
      \begin{smallmatrix}
        u^3 & u \\
        u & u^3
      \end{smallmatrix}\right)$ joins $E$ and\/~$E^S$.
    Each one of the blocks $M$, $M^\delta$, and $M^{\delta^2}$ is
    fixed by~$S$. No two of these three blocks have a point in common.
  \item %
    The block $I \coloneqq \left(
      \begin{smallmatrix}
        u^5 & 1 \\
        u & u^5
      \end{smallmatrix}\right)$
    joins $\p E$ with~$\p E^{S\delta}$, the block~$I^\delta = \left(
      \begin{smallmatrix}
        u^6 & 1 \\
        u^4 & u^6
      \end{smallmatrix}\right)$
    joins~$\p E^\delta$ with~$\p E^{S\delta^2}$, and the
    block~$I^{\delta^2} = \left(
      \begin{smallmatrix}
        u^3 & 1 \\
        u^2 & u^3
      \end{smallmatrix}\right)$
    joins~$\p E^{\delta^2}$ with~$\p E^{S}$.
  \item %
    The blocks~$I^\delta$ and~$I$ have the point~$\p F$ in common,
    where $F \coloneqq (I^\delta I)^3 %
    = \left(
      \begin{smallmatrix}
        u^6 & u^6 \\
        u^4 & u^2
      \end{smallmatrix}\right)$. %
  \item The points incident with~$I$ are $\p E$, $\p E^{S\delta}$,
    $\p F$, and $\p F^\delta$, where $F^\delta = \left(
      \begin{smallmatrix}
        u^3 & u^3 \\
        u^2 & u
      \end{smallmatrix}\right)$.
  \item The blocks joining~$\p D$ with points on~$I$ are $T$, $T^S$,
    $D\vee F = \left(
      \begin{smallmatrix}
        u^2 & u \\
        u^4 & u^2
      \end{smallmatrix}\right)$, and
    $D \vee F^\delta = (D \vee F^\delta)^\delta = \left(
      \begin{smallmatrix}
        u & u^4 \\
        u^2 & u
      \end{smallmatrix}\right)$.
  \item\label{oNan}%
    An O'Nan configuration is formed by the six points $\p D$,
    $\p E$, $\p {E^\delta}$, $\p E^{S\delta}$, $\p E^{S\delta^2}$,
    $\p F$, together with the four blocks $T$, $T^S$, $I$, $I^\delta$.
    That configuration is the dual of the complete graph~$K_4$ on $4$
    vertices.

    \begin{figure}
    \centering
    \begin{tikzpicture}[%
      every node/.append style={circle, 
        %draw=black, fill=green, %
        inner sep=1pt, %
        minimum size=0pt}]%

      \def\cw{white}%
      \def\cb{black}%{blue}
      \def\cr{black}%{red}
      \def\co{lightgray}%{orange}
      \def\cg{black}%{green}
      \def\cd{black}%color for D
      
      \def\r{.45}%
      \def\s{.75}%
      \coordinate (D) at (0,0) {};%
      \coordinate (E) at (-2,-6) {};%
      \coordinate (Ed) at ($(D)!\s!(E)$) {};%
      \coordinate (Edd) at ($(D)!\r!(E)$) {};%
      \coordinate (ES) at (2,-6) {};%
      \coordinate (ESd) at ($(D)!\s!(ES)$) {};%
      \coordinate (ESdd) at ($(D)!\r!(ES)$) {};%
      \coordinate (F) at (intersection of E--ESd and Ed--ESdd) {};%
      \coordinate (Fd) at (intersection of Ed--ESdd and Edd--ES) {};%
      \coordinate (Fdd) at (intersection of Edd--ES and E--ESd) {};%
      \coordinate (FS) at (intersection of ES--Ed and ESd--Edd) {};%
      \coordinate (FSd) at (intersection of ESd--Edd and ESdd--E) {};%
      \coordinate (FSdd) at (intersection of ESdd--E and ES--Ed) {};%
      \coordinate (G) at (0,-2.1) {};%
      \coordinate (Gd) at ($(D)!.33!(G)$) {};%
      \coordinate (Gdd) at ($(D)!.66!(G)$) {};%

      \draw[densely dashed] (D) -- (G) ; %
      \draw[densely dashed] (E) -- (ES) ; %
      \draw[densely dashed] (Ed) -- (ESd) ; %
      \draw[densely dashed] (Edd) -- (ESdd) ; %
      \draw[densely dashed] (F) -- (FS) ; %
      \draw[densely dashed] (Fd) -- (FSd) ; %
      \draw[densely dashed] (Fdd) -- (FSdd) ; %
      \draw (D) -- (E) ; %
      \draw (D) -- (ES) ; %
      \draw (E) -- (ESd) -- (F) -- (Fdd) -- cycle ; %
      \draw (Ed) -- (ESdd) -- (Fd) -- (F) -- cycle ; %
      \draw (Edd) -- (ES) -- (Fdd) -- (Fd) -- cycle ; %
      \draw[\cr] (ES) -- (Ed) -- (FS) -- (FSdd) -- cycle ; %
      \draw[\cr] (ESd) -- (Edd) -- (FSd) -- (FS) -- cycle ; %
      \draw[\cr] (ESdd) -- (E) -- (FSdd) -- (FSd) -- cycle ; %
      
      \node[draw=black,fill=\cd] at (D) {};
      \node[draw=black,fill=\cb] at (E) {};
      \node[draw=black,fill=\cb] at (Ed) {};
      \node[draw=black,fill=\cb] at (Edd) {};
      \node[draw=black,fill=\cg] at (F) {};
      \node[draw=black,fill=\cg] at (Fd) {};
      \node[draw=black,fill=\cg] at (Fdd) {};
      \node[draw=black,fill=\cr] at (ES) {};
      \node[draw=black,fill=\cr] at (ESd) {};
      \node[draw=black,fill=\cr] at (ESdd) {};
      \node[draw=black,fill=\co] at (FS) {};
      \node[draw=black,fill=\co] at (FSd) {};
      \node[draw=black,fill=\co] at (FSdd) {};
      \node[draw=black,fill=\cw] at (G) {};
      \node[draw=black,fill=\cw] at (Gd) {};
      \node[draw=black,fill=\cw] at (Gdd) {};

      \node[anchor=south] at (D) {$\p D$};%
      \node[anchor=north] at (E) {$\p E$};%
      \node[anchor=south east] at (Ed) {$\p E^\delta$};%
      \node[anchor=south east] at (Edd) {$\p E^{\delta^2}$};%
      \node[anchor=north] at (ES) {$\p E^S$};%
      \node[anchor=south west] at (ESd) {$\p E^{S\delta}$};%
      \node[anchor=south west] at (ESdd) {$\p E^{S\delta^2}$};%
      \node[anchor=south] at (F) {$\p F$};%
      \node[anchor=east] at (Fd) {$\p F^\delta$};%
      \node[anchor=north] at (Fdd) {$\p F^{\delta^2}$};%
      \node[anchor=south] at (FS) {$\p F^S$};%
      \node[anchor=west] at (FSd) {$\p F^{S\delta}$};%
      \node[anchor=north] at (FSdd) {$\p F^{S\delta^2}$};%
      \node[anchor=north] at (G) {$\p G$};%

      \node[fill=white] at ($(D)!.2!(E)$) {$T$};%
      \node[fill=white,rectangle,shift={(3pt,1pt)}] at ($(D)!.2!(ES)$) {$T^S$};%
      \node[fill=white] at ($(E)!.5!(F)$) {$I$};%
      \node[fill=white] at ($(Ed)!.5!(F)$) {$I^\delta$};%
      \node[fill=white,rectangle] at ($(ES)!.5!(FS)$) {$I^S$};%
      \node[fill=white,yshift=3pt,rectangle] at ($(ESd)!.5!(FS)$) {$I^{S\delta}$};%
    \end{tikzpicture}
    \caption{A super O'Nan configuration (black points), and its mirror image}
    \label{fig:twoSuperONanConfigurations}
  \end{figure}
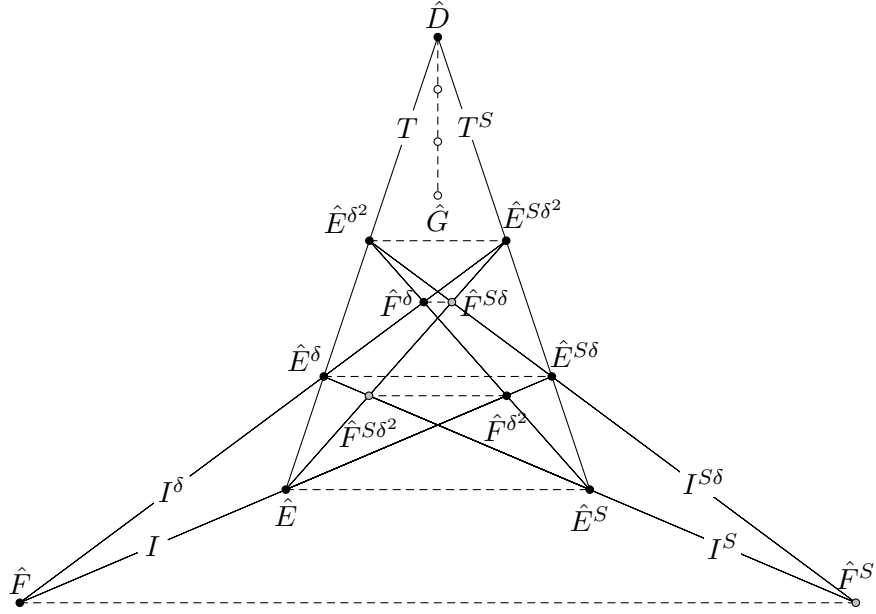
  \item\label{dualK5}%
    The ten points $\p D$, %
    $\p E$, $\p E^\delta$, $\p E^{\delta^2}$, %
    $\p E^S$, $\p E^{S\delta}$, $\p E^{S\delta^2}$, %
    $\p F$, $\p F^\delta$, $\p F^{\delta^2}$, %
    together with the five blocks\break%
    $T$, $T^S$, $I$, $I^\delta$, $I^{\delta^2}$ form a configuration
    which is the dual of the complete graph~$K_5$ on $5$ vertices. See
    Figure~\ref{fig:twoSuperONanConfigurations}. %

    This is one of the ``super O'Nan configurations'' found by
    Brouwer~\cite{MR655065}, see~\cite[Sect.\,4]{MR4264172}.
  \item Applying~$S$, we obtain a second super O'Nan configuration,
    with points $\p D$, %
    $\p E$, $\p E^\delta$, $\p E^{\delta^2}$, %
    $\p E^S$, $\p E^{S\delta}$, $\p E^{S\delta^2}$, %
    $\p F^S$, %
    $\p F^{S\delta}$, %
    $\p F^{S\delta^2}$, %
    and blocks %
    $T$, $T^S$, $I^S = \left(
      \begin{smallmatrix}
        u^5 & u \\
        1 & u^5
      \end{smallmatrix}\right)$, %
    $I^{S\delta} = \left(
      \begin{smallmatrix}
        u^6 & u^4 \\
        1 & u^6
      \end{smallmatrix}\right)$, %
    $I^{S\delta^2} = \left(
      \begin{smallmatrix}
        u^3 & u^2 \\
        1 & u^3
      \end{smallmatrix}\right)$,
    where
        $F^S = \left(
      \begin{smallmatrix}
        u^2 & u^4 \\
        u^6 & u^6
      \end{smallmatrix}\right)$, %
    $F^{S\delta} = \left(
      \begin{smallmatrix}
        u & u^2 \\
        u^3 & u^3
      \end{smallmatrix}\right)$, and %
    $F^{S\delta^2} = \left(
      \begin{smallmatrix}
        u^4 & u \\
        u^5 & u^5
      \end{smallmatrix}\right)$. %

  \begin{figure}
    \centering
    \begin{tikzpicture}[%
      scale=1.4, 
      every node/.append style={circle, 
        draw=black,
        fill=black, %
        inner sep=0pt, %
        minimum size=24pt}%
      ]%

      \def\mR{1.2}%
      \coordinate (T)    at (0,0) {};%
      \coordinate (TS0)  at (-2*\mR,0) {};%
      \coordinate (TS1)  at (\mR,{\mR*sqrt(3)}) {};%
      \coordinate (TS2)  at (\mR,{-\mR*sqrt(3)}) {};%
      \coordinate (I)    at (0,2) {};%
      \coordinate (Id)   at ({-sqrt(3)},-1) {};%
      \coordinate (Idd)  at ({sqrt(3)},-1) {};%
      \coordinate (IS)   at (0,-2) {};%
      \coordinate (ISd)  at ({sqrt(3)},1) {};%
      \coordinate (ISdd) at ({-sqrt(3)},1) {};%

      \coordinate (E)    at ($(I)!.5!(ISdd)$)  {};%
      \coordinate (Ed)   at ($(Id)!.5!(IS)$)   {};%
      \coordinate (Edd)  at ($(Idd)!.5!(ISd)$) {};%
      \coordinate (ES)   at ($(IS)!.5!(Idd)$)  {};%
      \coordinate (ESd)  at ($(ISd)!.5!(I)$)   {};%
      \coordinate (ESdd) at ($(ISdd)!.5!(Id)$) {};%

      \coordinate (F)    at ($(Idd)!.5!(I)$) {};%
      \coordinate (Fd)   at ($(I)!.5!(Id)$)  {};%
      \coordinate (Fdd)  at ($(Id)!.5!(Idd)$)   {};%
      \coordinate (FS)   at ($(ISdd)!.5!(IS)$) {};%
      \coordinate (FSd)  at ($(IS)!.5!(ISd)$)  {};%
      \coordinate (FSdd) at ($(ISd)!.5!(ISdd)$)   {};%

      \draw[lightgray,densely dotted,ultra thick] (T) circle (2*\mR) ;%
      \draw[ultra thick] (I) to (Id) to (Idd) to cycle ;%
      \draw[ultra thick] (IS) to (ISd) to (ISdd) to cycle ;%
      \draw[ultra thick] (I) to (Id) to (Idd) to cycle ;%
      \draw[ultra thick] (IS) to (ISd) to (ISdd) to cycle ;%

      \draw[ultra thick] (I) to (ISdd) to (Id) to (IS) to (Idd) to (ISd) to cycle ;%
      \draw[ultra thick] (T) to (E) ;%
      \draw[ultra thick] (T) to (Ed) ;%
      \draw[ultra thick] (T) to (Edd) ;%
      \draw[ultra thick] (TS2) to (ES) ;%
      \draw[ultra thick] (TS1) to (ESd) ;%
      \draw[ultra thick] (TS0) to (ESdd) ;%

      \node[fill=white] at (I)    {$I^{}$} ;%
      \node[fill=white] at (Id)   {$I^{\delta}$} ;%
      \node[fill=white] at (Idd)  {$I^{\delta^2}$} ;%
      \node[fill=white] at (IS)   {$I^{S}$} ;%
      \node[fill=white] at (ISd)  {$I^{S\delta}$} ;%
      \node[fill=white] at (ISdd) {$I^{S\delta^2}$} ;%

      \node[fill=white] at (T)    {$T$};%
      \node[fill=white,densely dotted] at (TS0)    {$T^S$};%
      \node[fill=white,densely dotted] at (TS1)    {$T^S$};%
      \node[fill=white,densely dotted] at (TS2)    {$T^S$};%
      
      \node at (E)    {\color{white}$\p E^{}  $};%
      \node at (Ed)   {\color{white}$\p E^{\delta} $};%
      \node at (Edd)  {\color{white}$\p E^{\delta^2}$};%
      \node at (ES)   {\color{white}$\p E^{S} $};%
      \node at (ESd)  {\color{white}$\p E^{S\delta}$};%
      \node at (ESdd) {\color{white}$\p E^{S\delta^2}$};%

      \node at (F)    {\color{white}$\p F^{}  $};%
      \node at (Fd)   {\color{white}$\p F^{\delta} $};%
      \node at (Fdd)  {\color{white}$\p F^{\delta^2}$};%
      \node at (FS)   {\color{white}$\p F^{S} $};%
      \node at (FSd)  {\color{white}$\p F^{S\delta}$};%
      \node at (FSdd) {\color{white}$\p F^{S\delta^2}$};%
    \end{tikzpicture}
    \caption{The incidence graph of the configuration $\cK$: the
      block~$T^S$ shows up in three places that have to be
      identified, %
      and the point $\p D$ is not shown.  (A spatial model would be
      better, but we do not have a 3D printer at hand.)}
    \label{fig:incidence graph}
  \end{figure}
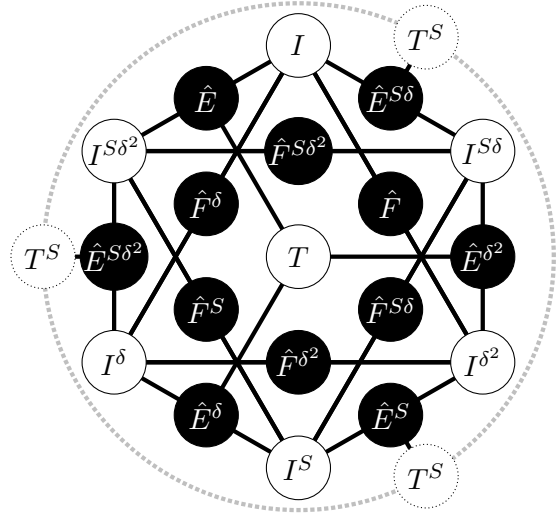

\item\label{configurationK}%
    The blocks $I^{S\delta}$ and $I^{S\delta^2}$ meet~$I$ in the
    points~$\p E^{S\delta}$ and $\p E$, respectively. %
    The blocks $I^S$ and $I$ have no point in common (as $I^SI$ has
    trace~$u^6$). %
    For any point $\p X$ not incident with the block~$S$, the block
    $B_X$ joining $\p X$ with $\p X^S$ is fixed by~$S$. The point rows
    of the six blocks $B_X$ with
    $X\in\{E,E^\delta,E^{\delta^2},F,F^\delta,F^{\delta^2}\}$ form a
    partition of the set of points not incident with~$S$.
    \SaveEnumi
  \end{enumerate}
  The union of the two super O'Nan configurations is a
  configuration~$\cK$ with $13$ points and~$8$ blocks, apart from the
  special point $\p D$, there are~$6$ points (of type~``$F$'')
  incident with $2$ blocks, and~$6$ points (of type~``$E$'') incident
  with~$3$ blocks.  The blocks come in two types: the blocks of
  type~``$I$'' are incident with $2$ points of type~``$E$'' and $2$ of
  type~``$F$'', while the blocks~$T$ and~$T^S$ are incident with $3$
  points of type~``$E$'' and the special point~$\p D$. See
  Figure~\ref{fig:twoSuperONanConfigurations} and
  Figure~\ref{fig:incidence graph}.
  \begin{enumerate}
    \RecallEnumi
  \item%
    The configuration~$\cK$ consists of the orbits of $\p D$, $\p E$,
    $\p F$, $T$, $I$, respectively, under the cyclic group
    $\langle\delta,S\rangle$ of order~$6$.
  \end{enumerate}
\end{exas}

\begin{ndef}[Open Problems]
  Is the Ree-Tits unital of order~$3$, among all unitals of order~$3$,
  characterized by the existence of configurations as
  in~\ref{explicit}\ref{dualK5} or
  in~\ref{explicit}\ref{configurationK} ?

  Consider a positive integer~$\ell$, the Ree-Tits
  unital~$\RT{3^{2\ell+1}}$ of order~$3^{2\ell+1}$ (see~\ref{def:RT}
  below), and an integer $n$ with $5<n\le 3^{2\ell+1}+2$.  Does there
  exist a configuration in~$\RT{3^{2\ell+1}}$ which is the dual of the
  complete graph on $n$ vertices?
\end{ndef}

\goodbreak%
\section{Ree groups}

The Ree-Tits unitals are incidence geometries that are closely related
to the Ree groups, see~\ref{def:RT} below. We collect some information
about the Ree groups first.

In order to construct a Ree group $\Ree[\theta]\KK$ (of type $\tGtwo$)
one needs a field $\KK$ of characteristic~$3$ and a \emph{Tits
  endomorphism}~$\theta$, i.e., an endomorphism $\theta$ of $\KK$ with
$\theta^2=\phi$, where $\phi\colon x\mapsto x^3$ denotes the Frobenius
endomorphism. %
We have $\theta = \id = \phi$ if, and only if, the field~$\KK$ is the
prime field~$\FF_3$. %
In general, existence and uniqueness of Tits endomorphisms depend on the
structure of~$\KK$. For a finite field $\KK$ of order~$3^n$, a Tits
endomorphism exists precisely if~$n$ is odd; in fact one has
$\theta=\phi^{(n+1)/2}$.

\begin{ndef}[Ree groups]
  We take the construction from~\cite{MR2795696}, following
  Tits~\cite{MR1611778} and~\cite{MR2730410}. %
  Write $\elt abc$ for the column with entries $a, b, c \in\KK$, and
  define a group operation~$*$ on the set of columns by
  \[
    \Elt abc*\Elt xyz := \EltAl%
    {{}{a+x}{}}
    {{}{b+y}{+ax^\theta}}
    {{ay-bx+}{c+z}{-ax^{\theta+1}}} \,.
  \]
  We denote\footnote{\ Other authors denote that ``root group'' by
    $\RUold[\theta]\KK$.} %
  this group by $\Xi \coloneqq \RU[\theta]{\KK}$.  If $\KK$
  is finite of order~$q$ we write~$\RU{q}:=\RU[\theta]{\KK}$.

  The following transformation~$\omega$ of the set of non-trivial
  elements of $\Xi$ is taken from~\cite{MR1611778} (with a correction
  in the definition of~$N$, cp.~\cite{MR2730410}):
  \[
  \omega\colon \Elt abc \mapsto \dfrac{-1}{N(a,b,c)}
  \Elt[1ex]{{a^\theta b^\theta - c^\theta + ab^2 + bc - a^{2\theta+3}}}
  {{ a^2 b - ac + b^\theta - a^{\theta+3}}} {{ c}}
  \]
  where $N(a,b,c) \coloneqq N(\elt abc) \coloneqq %
  -ac^\theta + a^{\theta+1}b^{\theta} - a^{\theta+3}b - a^2b^2 +
  b^{\theta+1} + c^2 - a^{2\theta+4}$. %
  (In~\cite[5.3]{MR1611778} one has to correct a misprint: replace the
  summand $a^{\theta+1}b$ by $a^{\theta+1}b^{\theta}$.) %
  
  Adding a new symbol~$\infty$ to the set $\Xi$, we extend~$\omega$ to
  a transformation of the set $P:= \Xi \cup \{\infty\}$; the elements
  $\infty$ and $o \coloneqq \elt000$ are swapped by~$\omega$. We note
  that~$\omega$ is an involution on~$P$ (see~\cite[p.\,16]{MR2730410}), but
  not an automorphism of~$\Xi$.

  We consider $\Xi$ as a group of permutations of~$P$, acting by
  multiplication from the right on itself (and fixing~$\omega$). %
  For each $s\in\KK^\times$, define $\eta_s \in \Aut{\Xi}$ by
  \[
    \Elt abc^{\eta_s} \coloneqq \EltAl{{s\ }a{}}{{s^{\theta+1}}b{}}{{s^{\theta+2}}c{}} .
  \]
  These maps will also be considered as permutations of~$P$,
  fixing~$\infty$. %
  Note that $N(\xi^{\eta_s}) = s^{2\theta+4}N(\xi)$. %
  Clearly, the set $\Eta \coloneqq \set{\eta_s}{s\in\KK^\times}$ is a
  group isomorphic to~$\KK^\times$. %
  Let $\KK^\dagger$ be the subgroup of~$\KK^\times$ generated by
  $N(\KK^3)\smallsetminus\{0\}$, and put
  $\Eta^\dagger \coloneqq \set{\eta_s}{s\in \KK^\dagger}$. %
  If $\KK$ is finite then $\Eta^\dagger = \Eta$ but
  $\Eta/\Eta^\dagger$ may be infinite in general
  (see~\cite[\S\,7]{MR2730410}). %
  Note that $-1 = N(0,1,1) \in \KK^\dagger$ yields
  $\eta_{-1}\in\Eta^\dagger$, in any case. From $N(0,0,c)=c^2$ we
  obtain that $\KK^\dagger$ contains all squares in~$\KK^\times$.
  
  The \emph{Ree group} $\Ree[\theta]\KK$ is the group of bijections
  of~$P$ generated by the subset
  $\{\omega\}\cup\Xi = \{\omega\} \cup \RU[\theta]\KK$ of $\Sym{P}$.
  We denote the Ree group by $\Ree{q}$ if $\KK$ is finite of
  order~$q$.
\end{ndef}

\begin{rema}\label{rem:computations}
  Some care is needed because we use elements of~$\Xi$ in different
  roles: as points in $P = \{\infty\}\cup\Xi$, and as permutations
  of~$P$. To wit, we chose the point
  $o \coloneqq \elt000 \in P\smallsetminus\{\infty\}$; then the point
  $\xi\in\Xi$ is $o^\xi$, where $\xi$ is interpreted as a permutation.

  We write $\xi^\omega$ and $\xi^{\eta_s}$ for the application
  of~$\omega$ and of~$\eta_s\in\Eta$, respectively, to $\xi\in\Xi$.
  For $\alpha,\beta\in\Xi$ and $\eta_s\in\Eta$, the product
  $\alpha\eta_s\omega\beta$
  in the group $\Ree[\theta]\KK$ then maps the point $\xi = o^\xi$ to
  $\bigl((\xi\alpha)^{\eta_s}\bigr)^\omega\beta =
  \bigl(\bigl((o^{\xi\alpha})^{\eta_s}\bigr)^\omega\bigr)^\beta =
  \bigl((o^{\eta_s^{-1}\xi\alpha\eta_s})^\omega\bigr)^\beta$; we may
  interpret $(\xi\alpha)^{\eta_s} = \xi^{\eta_s} \alpha^{\eta_s}$ as
  conjugation in the semidirect product $\Eta\Xi$.

  Note, however, that~$\omega$ is not an automorphism of~$\Xi$; this
  is where we need some care (and some parentheses).
\end{rema}

We quote the pertinent results from~\cite[Thm.\,1.1]{MR2730410}:

\begin{theo}\label{pertinentResults}
  \begin{enumerate}
  \item\label{pertinentResults1}%
    The group $\Ree[\theta]\KK$ acts doubly transitively on~$P$.
  \item\label{pertinentResults2}%
    The group $\RU[\theta]\KK$ is a normal subgroup of the
    stabilizer~$\Ree[\theta]\KK_\infty$ of\/~$\infty$.
  \item\label{pertinentResults3}%
    The stabilizer is a semidirect product\/
    $\Ree[\theta]\KK_\infty = \Eta^\dagger\RU[\theta]\KK$.
  \item\label{pertinentResults4}%
    Conjugation by\/~$\omega$ inverts each element of\/~$\Eta$. %
  \end{enumerate}
\end{theo}

\begin{rems}
  For the finite case (where $\KK$ has order~$q = 3^{2n+1}$), the Ree
  groups are discussed in~\cite{MR0138680} and in~\cite{MR0193136},
  see also~\cite[7.7.10]{MR1725957}.  One finds
  $|\Ree{q}|=(q^3+1)q^3(q-1)$; and\/ $\RU{q}$ is a Sylow $3$-subgroup
  of\/~$\Ree{q}$.
\end{rems}

\begin{ndef}[The nilpotent radical of the stabilizer]
  \label{stabilizer}
  For arbitrary $a,b,c\in\KK$, we have %
  \[
    {\Elt abc}^{-1} =
    \EltAl
    {{}{-a}{}}
    {{a^{\theta+1}}{-b}{}}{{}{-c}{}} , %
    \text{ and } %
    {\Elt abc}^{3} =
    \Elt{0}{0}{a^{\theta+2}}.
  \]
  In particular, the group $\Xi$ has exponent~$3^2$.  The set
  $\Lambda \coloneqq \C[\Xi]{\eta_{-1}} = \set{\elt0b0}{b\in\KK}$
  forms a subgroup of~$\Xi$, isomorphic to~$\KK$.  %

  The commutator
  $[\alpha,\xi] \coloneqq (\xi\alpha)^{-1}(\alpha\xi) =
  \alpha^{-1}\xi^{-1}\alpha\xi$ of $\alpha = \elt abc$ and
  $\xi = \elt xyz$ in~$\Xi$ is %
  \[
    \left[\Elt abc,\Elt xyz\right] = %
    \Elt0{x^\theta a-a^\theta x}{ay-bx+(x-a)(a^\theta x+x^\theta a)} .
  \]
  In particular, we have
  \[
      \left[\Elt 100,\Elt0y0\right] = %
    \Elt00y 
    \quad\text{ and }\quad%
      \left[\Elt 100,\Elt x00\right] = %
    \Elt0{x^\theta - x}{(x-1)(x+x^\theta)} .
  \]  
  In any case, the center of~$\Xi$ is
  $\Zeta \coloneqq \set{\elt00c}{c\in\KK}$, and we see that~$\Xi$ is nilpotent.

  If $|\KK|=3$ then $\Zeta = \Xi'$, and $\Xi$ is nilpotent of
  class~$2$. If $|\KK| > 3$ then $\Zeta$ is properly contained
  in~$\Xi'$, and the nilpotency class is~$3$. 
\end{ndef}

For the sake of completeness, we determine commutator groups:

\begin{prop}\label{commutators}
  Let\/ $\theta$ be a Tits endomorphism of a field\/~$\KK$ with
  $\Char\KK=3$. %
  \begin{enumerate}
  \item\label{order3commutators}%
    If\/ $|\KK|=3$ then
    $\RU[\theta]\KK' = \set{\elt00c}{c\in\KK}$, and the commutator
    group of the stabilizer\/ $\Ree[\theta]\KK_\infty$ equals
    $\set{\elt a{-a^2}c}{a,c\in\KK}$.
  \item\label{orderGeneralCommutators}%
    If\/ $|\KK|>3$ then $\RU[\theta]\KK' = \set{\elt0bc}{b,c\in\KK}$,
    and\/~$\RU[\theta]\KK$ is the commutator group of\/
    $\Ree[\theta]\KK_\infty$.
  \end{enumerate}
\end{prop}
\begin{proof}
  We abbreviate $\Xi\coloneqq\RU[\theta]\KK$ and $R
  \coloneqq\Ree[\theta]\KK$. %
  In any case, the group $(R_\infty)'$ contains 
  \[
    \begin{array}{rcccc}
      \left[\Elt abc,\eta_{-1}\right] %
      &=& \Elt abc^{-1} %
          * \Elt abc^{\eta_{-1}} %
      &=& \EltAl%
          {{}{-a}{}}%
          {{a^{\theta+1}}{-b}{}}%
          {{}{-c}{}} %
          *\Elt{-a}b{-c} %
      \\[4ex]
      &=& \EltAl%
          {{}{a}{}}%
          {{a^{\theta+1}-b}{+b}{+a^{\theta+1}}}%
          {{ab+a^{\theta+2}}{+c}{+a^{\theta+2}}}
      &=&
          \Elt{a}{-a^{\theta+1}}{ab+c-a^{\theta+2}} \,.
    \end{array}
  \]
  Using 
  $Y\coloneqq\set{x^\theta-x}{x\in\KK}$,
  and elements of~$\Xi'$ as obtained in~\ref{stabilizer},
  we obtain that the sets %
  \(%
  J := \set{\elt x{-x^{\theta+1}}z}{x,z\in\KK} %
  \) %
  and %
  \(%
  \Lambda = \set{\elt0y0}{y\in Y} %
  \) %
  are contained in $(R_\infty)'$.    
  
  If $|\KK|=3$ then $\theta=\id$, and $Y=\{0\}$.  Then
  $\Eta \Lambda = \langle\eta_{-1}\rangle \Lambda$ is a commutative
  subgroup of~$R_\infty$, and $J = \set{\elt x{-x^2}z}{x,z\in\KK}$ is
  a normal subgroup in~$R_\infty$.  Now $\Eta\Lambda$ forms a
  complement to~$J$, and we obtain
  $(R_\infty)' = \set{\elt x{-x^2}z}{x,z\in\KK}$.  The rest of
  assertion~\ref{order3commutators} is the observation $\Xi'=\Zeta$
  made for $|\KK|=3$ in~\ref{stabilizer}.
  
  Now assume $|\KK|>3$, then $Y$ contains some element $y\ne0$.  In
  the endomorphism ring of the multiplicative group of~$\KK$, we
  compute $(\theta-1)(\theta+1) = \theta^2-1 = 2$. So $\KK^{\theta+1}$
  contains the set of squares. That set of squares
  generates the additive group of~$\KK$ because $\Char\KK\ne2$. %
  Since~$\Xi$ is normalized by~$\Eta$, its commutator group~$\Xi'$ is
  also invariant under~$\Eta$, and we find that~$\Xi'$ contains
  $\set{\elt0bc}{b\in B,c\in\KK}$, where $B$ is additively generated
  by $\set{s^2y}{s\in\KK}$. %
  Now $B = \KK y = \KK$, and $\Xi' = \set{\elt0bc}{b,c\in\KK}$
  follows.

  Finally, $\Xi'\cup\set{\elt a{-a^{\theta+1}}{-a^{\theta+2}}}{a\in\KK} \subseteq (R_\infty)'$
  implies $(R_\infty)' = \Xi$, as claimed.
\end{proof}

\begin{lemm}\label{NCxi}
  Let\/ $R \coloneqq \Ree[\theta]\KK$.  %
  For each non-trivial\/ $\zeta$ in the center of\/~$\Xi$, the
  centralizer $\C[R]\zeta$ equals\/~$\Xi$, and\/
  $R_\infty = \N[R]{\C[R]\zeta}$.
\end{lemm}
\begin{proof}
  Since the group $\Xi$ acts sharply transitive on
  $P\smallsetminus\{\infty\}$, we know that~$\infty$ is the only point
  in~$P$ fixed by~$\zeta$.  So the centralizer of~$\zeta$ is contained
  in~$R_\infty$.
  
  The group~$\Eta \cong \KK^\times$ acts semi-regularly on
  $\Xi \smallsetminus\Lambda$ because $2+\theta$ is an
  automorphism of~$\KK^\times$ (with inverse $2-\theta$).  So the
  centralizer $\C[R]\zeta$ is contained in~$\Xi$, and then equal
  to~$\Xi$.

  As $\C[R]\zeta = \Xi$ fixes exactly one point in~$P$, the normalizer
  $\N[R]{\Xi}$ also fixes that point, and then coincides
  with~$R_\infty$.
\end{proof}

\section{Involutions, and the Ree-Tits unitals}
\begin{lemm}\label{singleConjugacyClassOfInvolutions}
  The involutions in $\Ree[\theta]\KK_\infty$ form a
  single conjugacy class under~$\Xi$.
\end{lemm}
\begin{proof}
  Every element
  of~$\Ree[\theta]\KK_\infty = \Eta^\dagger\Xi \le \Eta\Xi$ is of the
  form $\alpha\xi$ with $\alpha\in\Eta$ and $\xi \in \Xi$. We have
  $\id = (\alpha\xi)^2$ precisely if
  $\xi^{-1} = \alpha\xi\alpha = \alpha^2(\alpha^{-1}\xi\alpha)$.  From
  $\alpha^{-1}\xi\alpha \in \Xi$ we then obtain $\alpha^2=\id$. As
  $\Xi$ does not contain any involution, we find $\alpha\ne\id$, and
  that~$\alpha$ equals the unique involution $\eta_{-1}\in\Eta$. %
  Now
  $\xi^{-4}(\alpha\xi)\xi^4 = \alpha(\alpha^{-1}\xi^{-4}\alpha)\xi^5 =
  \alpha\xi^9 = \alpha$ follows from $\alpha^{-1}\xi\alpha =
  \alpha\xi\alpha = \xi^{-1}$ and $\xi^9=\id$.
\end{proof}

\begin{coro}
  The involutions with fixed points form a single conjugacy class
  in $\Ree[\theta]\KK$. %
  \qed 
\end{coro}

\begin{ndef}[Explicit description of involutions]
  \label{invExplicit}
  For $\xi = \elt abc \in\Xi$, we evaluate the condition
  $\elt{-a}{a^{\theta+1}-b}{-c} = \xi^{-1} = \eta_{-1}\xi\eta_{-1} =
  \elt{-a}{b}{-c}$ and obtain $b = -a^{\theta+1}$, while $a,c\in\KK$
  are arbitrary.  So $\eta_{-1}\xi$ is an involution precisely if
  $\xi \in J \coloneqq \set{\elt a{-a^{\theta+1}}c}{a,c\in\KK}$. %
  From~\ref{commutators} we know that this set~$J$ forms a subgroup
  of~$\Xi$ if, and only if, the field~$\KK$ has order~$3$ (viz., the
  Tits endomorphism is trivial). %

  A point $\xi\in P\smallsetminus\{\infty\}$ is fixed by $\eta_{-1}$
  precisely if $\xi$ is in the
  centralizer~$\Lambda= \C[\Xi]{\eta_{-1}}$ of~$\eta_{-1}$.  So the
  block $\Fix{\eta_1}$ equals
  $\{\infty\}\cup\Lambda$. From~\ref{pertinentResults}\ref{pertinentResults4}
  we know $\eta_{-1}^\omega = \eta_{-1}^{-1} =
  \eta_{-1}$. Thus~$\omega$ leaves that block invariant. %
  Now $\Lambda$ is a subgroup of~$\Xi$ acting regularly on the point
  set~$\Lambda$ while fixing~$\infty$, and $\omega$ swaps~$\infty$
  with $\id\in\Lambda$. So the stabilizer of the block
  $\Fix{\eta_{-1}}$ acts two-transitively on that block. %

  For $\xi = \elt a{-a^{\theta+1}}c\in J$, the involution
  $\eta_{-1}\xi$ has
  \[
    \Fix{\eta_{-1}\xi} = \{\infty\} \cup
    \set{\Elt{-a}y{ay-c}}{y\in\KK} . %
  \]
  For any two distinct involutions with fixed points, it follows that
  their sets of fixed points are distinct; in fact the sets
  $\Fix{\eta_{-1}\xi}\smallsetminus\{\infty\}$ with $\xi\in J$ induce
  a partition of $P\smallsetminus\{\infty\}$.
\end{ndef}

\begin{defi}\label{def:RT}
  The \emph{Ree-Tits unital} is the incidence structure
  $\RT[\theta]\KK \coloneqq (P,\cB,\in)$, where~$\cB$ is the set of
  all non-empty fixed point sets of involutions in $\Ree[\theta]\KK$.
\end{defi}

Clearly, the group $\Ree[\theta]\KK$ acts by automorphisms
of~$\RT[\theta]\KK$. Since that action is doubly transitive on~$P$,
any two points are joined by a unique block in that incidence
structure, by the last remark in~\ref{invExplicit}.  If $\KK$ has
finite order~$q$ then there are exactly~$q^2$ many involutions in
$\Ree{q}_\infty$. In~$\RT{q}$ we thus have exactly~$q^2$ blocks
through any given point.  It follows that $\RT{q}$ is a finite unital
of order~$q$.

For each $B\in\cB$, there is exactly one involution
$\sigma_B\in \Ree[\theta]\KK$ such that $B = \Fix{\sigma_B}$;
see~\ref{invExplicit}.  Consequently, the stabilizer of~$B$
in~$\Ree[\theta]\KK$ is the centralizer
$\C[{\Ree[\theta]\KK}]{\sigma_B}$. %

\begin{exas}\label{exas:sigmaTau}
  We abbreviate $\sigma \coloneqq \eta_{-1}$,
  $\zeta \coloneqq\elt 001$, and $\tau \coloneqq \sigma\zeta$. Then
  $\sigma$, $\tau$, and $\tau^\sigma = \zeta\sigma = \sigma\zeta^2$
  are involutions. We find
  $\Fix\sigma = \{\infty\}\cup\Lambda =
  \{\infty\}\cup\set{\elt0y0}{y\in\KK}$,
  $\Fix\tau = \{\infty\}\cup\set{\elt0y{-1}}{y\in\KK}$, and
  $\Fix{\tau^\sigma} = \{\infty\}\cup\set{\elt0y1}{y\in\KK}$.
\end{exas}

\begin{exas}\label{exas:omegaEpsilon}
  As we use~$\omega$ now, we have to be careful in our computations
  (see~\ref{rem:computations}). %

  We will only use elements from $\Ree3$ in this example, where
  $\theta=\id$. Then
  \[
    \Elt abc * \Elt xyz = \EltAl{{}{a+x}{}} {{}{b+y}{+ax}} {{ay-bx+}{c+z}{-ax^{2}}} \,;
  \]
  the formulae for~$N$ and~$\omega$
  simplify to %
  \(
    N\left(\elt abc\right) = %
    -ac - a^2b^2 + b^{2} + c^2 - a^{2}  %
  \)
  and %
  \[
    \Elt abc^\omega = %
    \dfrac{1}{ac + a^2b^2 - b^{2} - c^2 + a^{2}}
    \Elt{ab - c + ab^2 + bc - a} { a^2 b - ac + b - a^2} { c} .
  \]
  We find
  \[
    \Fix\omega = \left\{\Elt011,\Elt01{-1},\Elt110,\Elt{-1}10\right\}
    \quad\text{ (in the case where $|\KK|=3$)}.
  \]
  Now let $\alpha\in\Xi$, and consider the conjugate
  $\iota \coloneqq (\omega\alpha)\omega(\omega\alpha)^{-1}$. For
  $\xi\in\Xi$, we have
  $\xi\in\Fix\iota \iff \xi^\omega\alpha \in \Fix\omega \iff \alpha
  \in (\xi^\omega)^{-1}\Fix\omega$. %
  In order to find the block joining
  $\varepsilon \coloneqq \zeta^{-1} =  \elt00{-1}$
  with~$\varepsilon' \coloneqq \varepsilon^{\sigma}\lambda$ for
  $\lambda \coloneqq \elt010$, we compute
  $(\varepsilon^\omega)^{-1} = (\elt{-1}01)^{-1} = \elt11{-1}$ and
  $\bigl((\varepsilon')^\omega\bigr)^{-1} = (\varepsilon')^{-1} =
  \elt0{-1}{-1}$. Then we obtain the intersection
  \[
    (\varepsilon^\omega)^{-1}\Fix\omega \cap
    \bigl((\varepsilon')^\omega\bigr)^{-1}\Fix\omega =
    \left\{\Elt{-1}01\right\} .
  \]
  With $\alpha = \elt{-1}01$ we obtain the joining block
  \[
    \Fix\iota = \Fix{(\omega\alpha)\omega(\omega\alpha)^{-1}} = %
  \bigl(\Fix\omega\alpha^{-1}\bigr)^\omega = %
  \left\{ \varepsilon, \varepsilon', \psi, \psi\lambda \right\} \,,
  \]
  where $\psi \coloneqq \elt{-1}{-1}1$ and $\psi\lambda = \elt{-1}00$.

  Clearly, the blocks $\Fix\iota$ and
  $\Fix\iota\lambda = \Fix{\iota^\lambda}$ share the point
  $\psi\lambda$, and $\psi\lambda^2$ is on both $\Fix{\iota^\lambda}$
  and $\Fix{\iota^{\lambda^2}}$.  So the ten points $\infty$,
  $\varepsilon$, $\varepsilon^{\lambda}$, $\varepsilon^{\lambda^2}$,
  $\varepsilon^{\sigma}$, $\varepsilon^{\sigma\lambda}$,
  $\varepsilon^{\sigma\lambda^2}$, $\psi$, $\psi^{\lambda}$,
  $\psi^{\lambda^2}$ and the five blocks $\Fix{\tau}$,
  $\Fix{\tau^\sigma}$, $\Fix{\iota}$, $\Fix{\iota^\lambda}$,
  $\Fix{\iota^{\lambda^2}}$ form a super O'Nan configuration~$\cS$
  in~$\RT3$. %

  In fact, this configuration was obtained by translating the one
  from \ref{explicit}\ref{dualK5} via a suitable isomorphism from
  $\sL28$ onto~$\Ree3$. See Section~\ref{sec:isoPsL28Ree3} below. 
\end{exas}

\begin{exas}\label{stringOfPearls}
  We return to the general case of $\Ree[\theta]\KK$ now, keeping
  notation from~\ref{exas:omegaEpsilon}. The block joining the
  points~$\varepsilon$ and~$\varepsilon'$ is
  still~$\Fix\iota = (\Fix\omega\alpha^{-1})^\omega$, where
  $\iota = (\omega\alpha)\omega(\omega\alpha)^{-1}$ with
  $\alpha = \elt{-1}01 \in\Xi$. If~$|\KK|>3$ then that fixed
  point set contains more points, of course.

  The group $\Lambda$ is abelian, and centralizes each one of the
  elements~$\sigma$, $\tau$, $\tau^\sigma$. For each $t\in\KK$, the
  automorphism $\lambda_t \coloneqq \elt0t0 \in\Lambda$ thus leaves invariant the
  blocks~$\Fix\sigma$, $\Fix\tau$, $\Fix{\tau^\sigma}$, and maps the
  super O'Nan configuration~$\cS$ to another such
  configuration~$\cS^\lambda$. We obtain $|\KK|/3$ super O'Nan
  configurations sharing the two blocks $\tau$, $\tau^\sigma$, and
  their common point~$\infty$.

  If $\KK$ is finite of order~$q$, the union of those $q/3$
  configurations has $3q+1$ points, and $q+2$ blocks.
\end{exas}

\begin{ndef}[Further intersections]%
  For $\mu\in\Lambda$ and~$\iota$ as in~\ref{exas:omegaEpsilon}, one
  may ask whether the blocks~$\Fix\iota$ and
  $\Fix\iota^\mu = \Fix{\iota^\mu}$ share a point if
  $\mu\notin\langle\lambda\rangle$.
  In order to answer that question, we determine
  $\Fix\iota$ in~$\RT[\theta]\KK$, for general~$\KK$.
  
  For $t\in\KK$, we write $\lambda_t \coloneqq \elt0t0 \in\Lambda$. %
  Recall that~$\Lambda$ is abelian, and that
  $\Fix\tau = {\{\infty\}\cup \varepsilon\Lambda} %
  = \{\infty\}\cup\Lambda\varepsilon$. %
  We apply~$\omega\varepsilon$ to
  $\{\infty, \varepsilon\lambda_1 \} \subseteq \Fix{\tau}$ and
  obtain %
  $\left\{\infty,\varepsilon\lambda_1\right\}^\omega\varepsilon =
  \left\{\varepsilon,\varepsilon'\right\} \subseteq \Fix\iota$ with
  $\varepsilon' = \varepsilon^\sigma\lambda_1$, as
  in~\ref{exas:omegaEpsilon}.  As joining blocks are unique in
  $\RT[\theta]\KK$, we infer $\tau^{\omega\varepsilon} = \iota$, and
  $\Fix\iota = \Fix\tau^{\omega}{\varepsilon} =
  \left(\{\infty\}\cup\Lambda\varepsilon\right)^{\omega}{\varepsilon}
  = \{\varepsilon\} \cup (\Lambda\varepsilon)^\omega\varepsilon$.  Now
  \[
    (\Lambda\varepsilon)^{\omega} =
    \set{(\lambda_x\varepsilon)^{\omega}}{\strut x\in\KK} =
    \set{\Elt0x{-1}^{\omega}}{x\in\KK} =
    \set{\frac{-1}{x^{\theta+1}+1}\Elt{1-x}{x^\theta}{-1}}{x\in\KK} .
  \]
  Consider $m\in\KK\smallsetminus\{0\}$ and
  $\mu = \lambda_m \in\Lambda$.  Assume that there exists
  $\xi \in \Fix\iota\cap\Fix\iota\mu$.  We have
  $\varepsilon\mu = \mu\varepsilon$ and thus
  $(\Lambda\varepsilon)^{\omega}\varepsilon\mu =
  (\Lambda\varepsilon)^{\omega}{\mu\varepsilon}$. %
  If $\xi \in \{\varepsilon,\varepsilon\mu\}$ then either
  $\mu^{-1}$ or~$\mu$ lies in
  $(\Lambda\varepsilon)^{\omega}$, which is impossible. %
  So both $\xi\varepsilon^{-1}$ and
  $\xi\mu\varepsilon^{-1}$ lie in
  $(\Lambda\varepsilon)^\omega$, and there exist
  $x,s\in\KK$ such that 
  \[
    \frac{-1}{x^{\theta+1}+1}\Elt{1-x}{x^\theta}{-1} =
    \frac{-1}{s^{\theta+1}+1}\Elt{1-s}{s^\theta}{-1}*\Elt0m0 =
    \frac{-1}{s^{\theta+1}+1}\Elt{1-s}{s^\theta-m(s^{\theta+1}+1)}{(1-s)m-1} .    
  \]
  These equations do not have any solutions with $1\in\{x,s\}$.  For
  $s\in\{0,-1\}$ we obtain the expected solutions (leading to
  $(\xi,\mu) \in \{(\psi,\lambda^2),(\psi\lambda,\lambda)\}$) and no
  others.

  Using GAP~\cite{GAP4-2022}, the author has checked that, at least for
  $|\KK|<3^{11} = 177\,147$, there are no solutions with
  $s\notin \FF_3$. So, in these small cases, the blocks $\Fix\iota$
  and $\Fix{\iota^\mu}$ do not share any point if
  $\mu\in\Lambda\smallsetminus\{\lambda_0,\lambda_1,\lambda_{-1}\}$.
\end{ndef}

\section{Translating from $\sL28$ into $\Ree3$}
\label{sec:isoPsL28Ree3}

We use the fact that there exists an isomorphism
$\gamma\colon \sL28 \to \Ree3$ (see~\cite[Sect.\,4.5.4]{MR2562037},
cp.~\cite[Thm\,1,\,9.]{MR2653666}), and use notation as
in~\ref{nota:PgL28} and in~\ref{exas:sigmaTau}. %
Under~$\gamma$, the Sylow $3$-subgroup
$\Delta = \langle\delta,A\rangle$ is mapped to a conjugate of
$\Xi$. Without loss of generality, we may thus assume
$\Delta^\gamma = \Xi$. Then
$\p A = \N[\sL28]{\langle A\rangle} = \N[\sL28]{\Delta}$ is mapped to
$\p A^\gamma = \N[\Ree3]\Xi = \Ree3_\infty$, and~$S^\gamma$ is a
conjugate of~$\sigma = \eta_{-1}$ under~$\Xi$;
see~\ref{singleConjugacyClassOfInvolutions}. Without loss of
generality, we may thus also assume $S^\gamma = \sigma$. %
As~$\delta$ generates the centralizer of~$S$
in~$\Delta$, we obtain $\delta^\gamma\in\C[\Xi]\sigma = \Lambda$.  So
$\delta^\gamma = \lambda_v \coloneqq \elt0v0$ with some
$v\in\{1,-1\}$.

Next, we note that the involutions $S$, $T$, $T^S$ are just those
involutions in~$\p A$ that commute with~$\delta$.  So $T^\gamma$ lies
in the centralizer $\C[R_\infty]{\delta^\gamma}$. This implies
$T^\gamma \in \sigma\Zeta$. %
The group $\Eta$ centralizes~$\sigma$, normalizes~$\Xi$, and acts
transitively on the set of non-trivial elements of~$\Zeta$. So we may
assume $T^\gamma = \tau$; then $(T^S)^\gamma = \tau^\sigma$. 

\goodbreak%
Our isomorphism~$\gamma$ now induces an isomorphism
$\check\gamma\colon \Ss \to \RT3$. Under that isomorphism, the point
$\p D$ goes to~$\infty$, the blocks $S$, $T$, $T^S$ go to~$\sigma$,
$\tau$, $\tau^\sigma$, respectively.  The point $\p E$ lies on the
block~$T$, so it goes to a point on the block $\Fix\tau$, that is, to
a point of the form $\elt0t{-1}$.  Applying the group
$\Lambda \le \C[R_\infty]{\langle\sigma,\tau\rangle}$, we achieve that
$\p E$ is mapped to $\varepsilon = \elt00{-1}$. We apply
$\langle S,\delta\rangle^\gamma = \langle\sigma,\lambda_1\rangle$ and
find the orbit
\[
  E^{\gamma\langle\sigma,\lambda_v\rangle} = %
  \left\{ %
    \Elt00{-1}, \Elt001, \Elt01{-1}, \Elt011, \Elt0{-1}{-1},
    \Elt0{-1}1 %
  \right\}
\]

%%=====================================================================
%%%%%%%%%%%%%%%%%%%%%%%%%%%%%%%%%%%%%%%%%%%%%%%%%%%%%%%%%%%%%%%%%%%%%%%
% \bibliographystyle{mybibstyle-noISSN-noISBN}%%%% %%%%% %%%%% %%%%% %%%%
% \bibliography{myBibliography}  %%%%% %%%%% %%%%% %%%%% %%%%% %%%%% %%%%
%   \enlargethispage{5mm}%
%%%%%%%%%%%%%%%%%%%%%%%%%%%%%%%%%%%%%%%%%%%%%%%%%%%%%%%%%%%%%%%%%%%%%%%
\providecommand{\noopsort}[1]{}\def\cprime{$'$}
  \def\polhk#1{\setbox0=\hbox{#1}{\ooalign{\hidewidth
  \lower1.5ex\hbox{`}\hidewidth\crcr\unhbox0}}}

%%%%%%%%%%%% ~~~~~~~~~~~~~~~~~~~~~~~~~~~~~~~~~~~~~~~~~~~~~~~~~~~~~~~~~~~
\medskip\vfill

\begin{small}
  \begin{minipage}[t]{0.3\linewidth}
    Markus J. Stroppel %<-------------------
  \end{minipage}
  \begin{minipage}[t]{0.6\linewidth}
    LExMath\\
    Fakult\"at 8\\
    Universit\"at Stuttgart\\
    70550 Stuttgart\\ %<-------------------
    stroppel@mathematik.uni-stuttgart.de %<-------------------
  \end{minipage}
\end{small}

%%%%%%%%%%%%~~~~~~~~~~~~~~~~~~~~~~~~~~~~~~~~~~~~~~~~~~~~~~~~~~~~~~~~~~~

\end{document}